# Extended Convex Hull-Based Distributed Optimal Energy Flow of Integrated Electricity-Gas Systems

Rong-Peng Liu[a], Wei Sun[b], Dali Zhou[c], Yunhe Hou[a,*]

[a] Department of Electrical and Electronic Engineering, the University of Hong Kong, Hong Kong 999077, China

[b] Department of Electrical and Computer Engineering, University of Central Florida, Orlando, FL. 32186, USA

[c] China National Heavy Duty Truck Group Co., Ltd., Jinan 250000, China

[*] Corresponding author: yhhou@eee.hku.hk

***Abstract***: Integrated electricity and gas systems are constructed to facilitate the gas-fired generation, and the distributed operation of these integrated systems have received much attention due to the increased emphasis on data security and privacy between different agencies. This paper proposes an extended convex hull based method to address optimal energy flow problems for the integrated electricity and gas systems in a distributed manner. First, a multi-block electricity-gas system model is constructed by dividing the whole system into $N$ blocks considering both physical and regional differences. This multi-block model is then convexified by replacing the nonconvex gas transmission equation with the extended convex hull-based constraints. The Jacobi-Proximal alternating direction method of multipliers algorithm is adopted to solve the convexified model and minimize its operation cost. Finally, the feasibility of the optimal solution for the convexified model is checked, and a sufficient condition is developed. If the sufficient condition is satisfied, the optimal solution for the original nonconvex problem can be recovered from that for the convexified problem. Simulation results demonstrate that the proposed method is tractable and effective in obtaining feasible optimal solutions for multi-block optimal energy flow problems.

***Keywords***: Alternating direction method of multipliers, convex relaxation, distributed optimization, integrated electricity and gas systems, optimal energy flow, quadratic programming

## Nomenclature

| | |
|---|---|
| *Abbreviations* | |
| ADMM | Alternating direction method of multipliers |
| CCP | Convex-concave procedure |
| ECH | Extended convex hull |
| IEGS | Integrated electricity and gas systems |
| J-ADMM | Jacobi-Proximal alternating direction method of multipliers |
| OEF | Optimal energy flow |
| PDAD | Primal-dual absolute difference |
| PWL | Piecewise linearization |
| QP | Quadratic programming |
| SOC | Second-order cone |
| *Sets* | |
| $\mathcal{D}_p^r/\mathcal{D}_g^{r'}$ | Set of power/gas loads in block $r/r'$ |
| $\mathcal{G}_p^r/\mathcal{G}_g^r/\mathcal{G}_v^{r'}$ | Set of coal-fired/gas-fired/virtual gas-fired generators in block $r/r/r'$ |
| $\mathcal{L}_p^r/\mathcal{L}_g^{r'}/\mathcal{C}_g^{r'}$ | Set of power transmission lines/gas passive pipelines/gas compressors (gas active pipelines) in block $r/r'/r'$ |
| $\mathcal{N}_p^r/\mathcal{N}_v^r/\mathcal{N}_g^{r'}$ | Set of power/virtual power/gas nodes in block $r/r/r'$ |
| $\mathcal{W}^{r'}$ | Set of gas wells in block $r'$ |
| $\mathcal{Y}_{\text{con}}/\mathcal{Y}_{\text{e-con}}$ | Convex hull/extended convex hull of set $\mathcal{Y}$ |
| $\Omega_p^r/\Omega_g^{r'}$ | Feasible region of power block *r*/gas block *r'* |
| *Parameters* | |
| $C_p(\cdot)$ | Quadratic cost function of coal-fired generator with quadratic term coefficient $c_1$, linear term coefficient $c_2$, and constant term $c_3$ |
| $C_w$ | Cost of natural gas |
| d | Penalty parameter |
| $G_g^{min}/G_g^{max}$ | Output limits of gas-fired generator |
| $G_i^{min}/G_i^{max}$ | Pressure square limits of gas node |
| $G_l^{min}/G_l^{max}$ | Transmission limits of gas passive pipeline |
| $G_w^{min}/G_w^{max}$ | Output limits of gas well |
| $N$ | Total number of blocks |
| $\mathbf{P}^r/\mathbf{E}^r$ | Positive semi-definite matrix/constant matrix |
| $P_d/G_d$ | Nodal loads of power/gas network |
| $P_g^{min}/P_g^{max}$ | Output limits of coal-fired generator |
| $P_l/G_c$ | Transmission limit of power transmission line/gas compressor |
| $W_l$ | Weymouth equation constant |
| $x_l$ | Reactance of power transmission line |
| $\alpha_c$ | Gas compressor constant |
| $\gamma$ | Damping parameter |
| $\theta_i^{min}/\theta_i^{max}$ | Phase angle limits of power node |
| $\chi_g$ | Electricity-gas conversion ratio |
| *Variables* | |
| $g_w$ | Output of gas well |
| $g_g^j$ | Output of virtual gas-fired generator |
| $p_g/g_g$ | Output of coal-fired/gas-fired generator |
| $p_l/g_l/g_c$ | Power/gas/gas flow through power transmission line/gas passive pipeline/gas compressor |
| $\boldsymbol{\delta}^+/\boldsymbol{\delta}^-$ | Slack variable (vector) |
| $\theta_i/\pi_i$ | Phase angle/pressure square of power/gas node |
| $\theta_i^j$ | Phase angle of virtual power node |
| $\lambda$ | Lagrangian multiplier |

## 1 Introduction

Traditional power generation relies heavily on coal and causes severe environmental problems. Seeking alternatives is imperative. Natural gas is a feasible solution due to its eco-friendly and low-cost properties. Thanks to the shale gas revolution and tax credit [1], natural gas has been one of the largest energy sources of energy for electricity generation, especially in the U.S. [2] and the U.K. [3] (38.4% and 40.6%, respectively, in 2019). Hereinafter, the word "gas" refers to both the underground natural gas and the natural gas extracted from the shale rock.

Bulk integrated electricity and gas system (IEGS) is constructed to facilitate gas transmission and gas-fired generation. As one of the most fundamental and critical problems for the IEGS, the optimal energy flow (OEF) problem has attracted many researchers' attention. Nevertheless, the IEGS brings new challenges caused by complex coupling relations, e.g., synergistic expansion, synchronous dispatch, and security related issues, in optimizing the energy flow. A large amount of research work was conducted to address these challenges. Specifically, reference [4] proposes a co-expansion planning framework to achieve a low-carbon economy for the IEGS. Reference [5] applies an expansion model to the IEGS considering bi-directional energy conversion. A practical OEF problem is proposed in [6], in which only a limited number of components can be adjusted. Reference [7] proposes a renewable penetrated IEGS model, in which power-to-gas conversion facilities are included to optimize the intraday operation strategy. In [8], security constraints are incorporated into the IEGS to ensure its feasibility under pre-defined security conditions. Reference [9] adopts a two-stage robust optimization model to enhance the resilience of the integrated electric-gas distribution systems against natural disasters. Please refer to [10] for a review of the integrated electricity and gas system coordination.

Previous work mainly focuses on centralized operation. In practice, power and gas networks usually belong to different companies. Even a connected power network may be divided into several blocks by region, with each block managed by an agency. From the perspective of security and privacy, it is risky and unrealistic to share all the information between different companies and agencies. *Distributed operation* is a promising solution and leads to distributed optimization problems. The alternating direction method of multipliers (ADMM) algorithm [11] is widely adopted to address power system optimization problems due to its promising performance. Reference [12] proposes a distributed algorithm based on ADMM to solve the nonconvex optimal power flow problem. Reference [13] applies the ADMM to handling wind power uncertainty in power systems and minimizing the operation cost. In [14], a multi-layer multi-agent model is proposed for optimal electric vehicle charging problems and is solved by the ADMM algorithm. Reference [15] proposes an energy management method to achieve the resilient and economic operation of the shipboard power system in a distributed manner. Reference [16] develops a method to adjust the penalty parameter automatically, which speeds up the ADMM convergence for the distributed optimal power flow problem. However, the following challenges still exist when addressing OEF problems for the IEGS by the ADMM algorithm: i) decoupling method for the IEGS; ii) convergence guarantee for the ADMM algorithm when the number of blocks is larger than two; iii) nonconvex constraints in the gas block.

Recently, some research work focuses on the above challenges. Specifically, based on the physical difference, reference [17] divides the IEGS into two blocks, i.e., power and gas blocks. The second-order cone (SOC) relaxation method [18] is adopted to relax nonconvex gas transmission equations, and the (relaxed) OEF problem is solved by the standard ADMM algorithm. In [19], during each iteration of the standard (two-block) ADMM algorithm, the local optimal solution for the nonconvex gas block is obtained by the convex-concave procedure (CCP) [20]. Reference [21] adopts the same two-block decoupling method and linearizes the nonconvex gas constraints by the piecewise linearization (PWL) method [22]. The OEF problem is solved by the tailored ADMM algorithm [23]. Differently, reference [24] divides the IEGS into $N$ ($N \geq 2$) blocks by region. The sequential cone programming method is leveraged to handle the nonconvex gas block, and the iterative ADMM is proposed to solve the OEF problem. Reference [25] proposes a two-block environmental-economic dispatch model. A two-layer distributed framework is adopted to address the OEF problem. In [26], a chance-constrained decentralized model is adopted to accommodate the renewable uncertainties and minimize the operation cost for the multi-region IEGS. In this paper, we attempt to address challenges i)-iii) in the following ways:

i) *Decoupling method:* Previous work employs the physical decoupling standard (in [17], [19], and [21]) and regional decoupling standard (in [24]) separately to build two or multi-block IEGS models. As aforementioned, even a connected power network may belong to different agencies, let along power and gas networks. The construction of the multi-block IEGS model considering both standards should be explored.

ii) *Distributed optimization:* Reference [27] states that the direct extension of the two-block ADMM algorithm to addressing multi-block distributed optimization problems does not necessarily converge. The Jacobi-Proximal ADMM (J-ADMM) algorithm [28], which is provably convergent when solving multi-block optimization problems, is more appropriate to solve distributed OEF problems for the multi-block IEGS.

iii) *Convex relaxation:* One of the convergence conditions for the distributed algorithm requires the convexity of each block [28]. However, the gas block model is nonconvex due to the gas transmission equation. Convex relaxation method proposed in [18] is a plausible solution and is adopted in [21] and [24]. The problem is that binary variables are introduced into the OEF problem, and its convergence cannot be guaranteed. Differently, references [19] and [25] adopt the CCP method to convexify the gas block model, while the bi-directional property of the gas flow can not be preserved anymore. In order to fill in this research gap, new convex relaxation methods for the distributed OEF problem are still worth exploring. Notice that the optimal solution for the relaxed problem may not be feasible for the original problem. The solution feasibility should be checked carefully.

Overall, this paper aims to address the distributed OEF problem for the *transmission-level* IEGS, and the main contributions are twofold:

1) Distributed optimal operation. A distributed framework is proposed to address the OEF problem for the multi-block IEGS. First, a decoupling method considering both physical and regional differences of the IEGS is developed to decouple the IEGS into $N$ ($N \geq 2$) blocks. Then, coupling relations between decoupled blocks are explored to formulate the distributed OEF problem for the multi-block IEGS. To the best of the authors' knowledge, this is the first work to construct a multi-block IEGS by leveraging both physical and regional differences. The J-ADMM algorithm, which allows parallel computing and can provide the unique optimal operation policy for the multi-block IEGS with guaranteed convergence (see Section 3.2 for details), is quite applicable and thus adopted to solve the distributed OEF problem.
2) Extended convex hull (ECH) based relaxation method. A convexified IEGS model is proposed by replacing nonconvex gas transmission equations with ECH-based constraints, which leads to the convexified OEF problem. Compared with [17] and [21], the proposed convexification method does not introduce binary variables and thus ensures the convergence of the distributed algorithm. Compared with [19] and [29], the bi-directional property of gas passive pipelines is preserved, which increases gas transmission flexibility. In addition, the feasibility of the solution for the convexified OEF problem is checked, and a sufficient condition is developed. The optimal solution for the original (nonconvex) OEF problem is recovered from that for the convexified OEF problem if the sufficient condition is satisfied. Simulation results demonstrate that this method is effective in obtaining feasible optimal solutions for radial gas networks.

The remainder of this paper is organized as follows. Section 2 introduces the decoupling method and the mathematical formulation of the multi-block IEGS. Section 3 presents the ECH, the J-ADMM algorithm, and the solution feasibility check and recovery method. Case studies are conducted in Section 4. Section 5 concludes this paper.

## 2 Problem Formulation

In this section, a decoupling method is proposed to construct the multi-block IEGS model, followed by the OEF formulation of power and gas blocks. Before that, we make the following assumptions and simplifications:

1) In order to reduce modelling complexity, this work focuses on the single-period OEF problem. In fact, the proposed distributed framework also applies to the multi-period OEF problem.
2) The IEGS is first divided into one power network and one gas network. Then, the power network is divided into multiple blocks. For the sake of simplicity, the gas network is no longer partitioned by assuming that it belongs to one gas company. Note that the proposed decoupling method for power networks can be directly extended to gas networks.
3) Power and gas load uncertainty is ignored, as this work aims to address the OEF problem in a distributed manner. Please refer to [13] (chance-constrained distributed optimal power flow with wind uncertainty), [30] (stochastic optimization-based distributed residential demand response with weather and consumer uncertainties), and [31] (robust optimization-based distributed reserve scheduling with multiple uncertainties) for ADMM-based methods utilized to address various uncertainties.
4) The OEF problem for the IEGS is feasible. This assumption is readily satisfied by setting load levels within a reasonable range. It follows the fact that a real-world IEGS which works normally can always satisfy its regular (power and gas) load demand and thus has at least one feasible solution for the OEF problem. In addition to our work, this assumption is also adopted in other research work, such as [32] and [33].

*2.1 Decoupling Method*

Based on the physical difference, the power network is decoupled from the gas network. Gas-fired generators are coupling components that connect these two networks, and the coupling relation is shown in the upper portion of Fig. 1. The dotted line denotes the gas supplied to the gas-fired generator $g$ from the gas node $j$. Gas-fired generators are considered to belong to power blocks. As is shown in the lower portion of Fig. 1, a virtual gas-fired generator $g'$ is added to the gas node $j$ to help decouple the power block from the gas block.

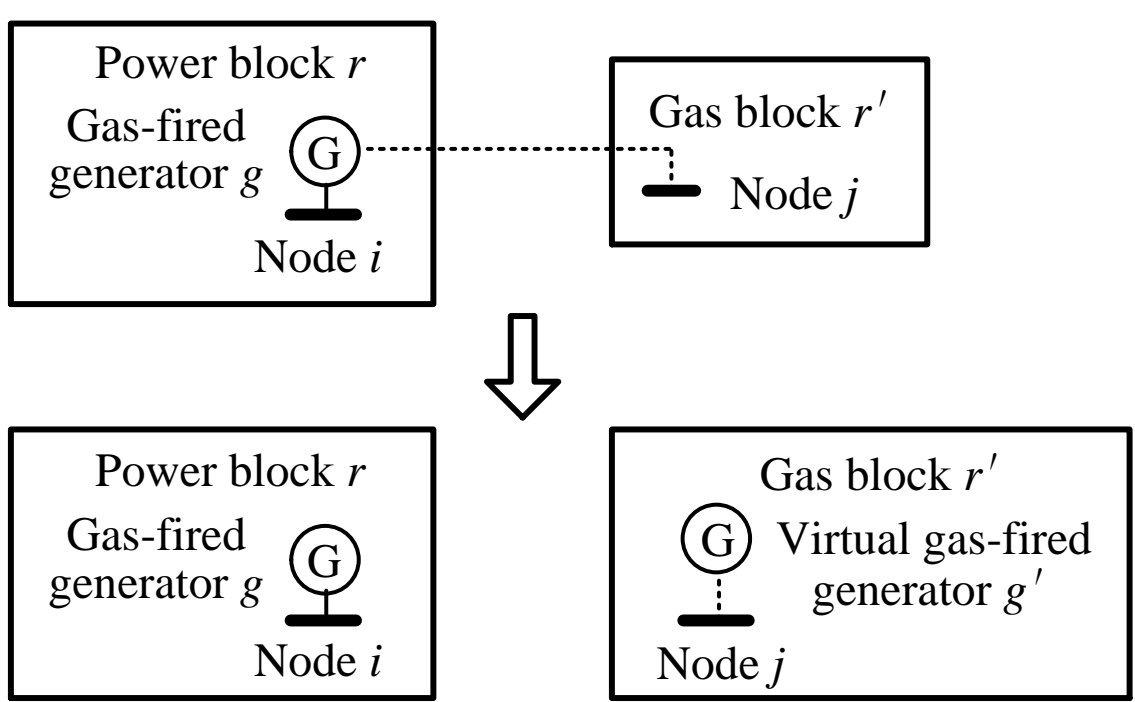

**Fig. 1.** Decoupling method between power block and gas block

Considering that power blocks in different regions usually belong to different agencies, the decoupled power network is further divided into several blocks based on the regional difference. Power transmission lines that connect different blocks are coupling components, and the coupling relation is shown in the upper portion of Fig. 2. Power blocks $r_1$ and $r_2$ are connected by a power transmission line $l$ (solid line). As is shown in the lower portion of Fig. 2, virtual power nodes $j'$ and $i'$ are added to power blocks $r_1$ and $r_2$, respectively. The original power transmission line $l$ is separated into two lines to connect the actual and the virtual nodes in these two blocks. By repeating this decoupling method, the power network is decoupled into multiple power blocks.

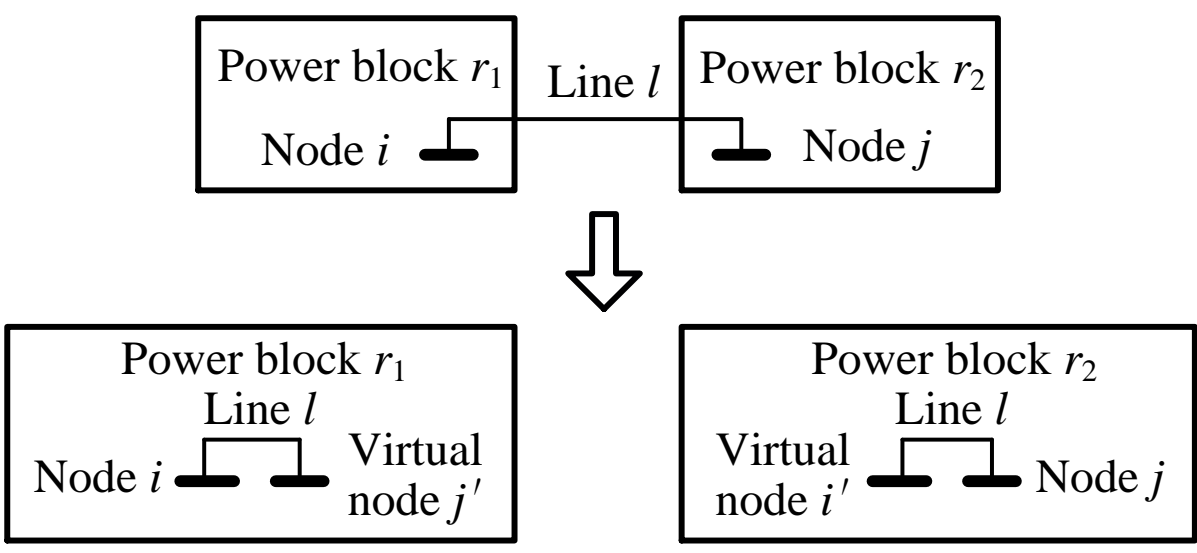

**Fig. 2.** Decoupling method between two power blocks

*2.2 Power Block Model*

Without loss of generality, we assume that the first $N$-1 blocks are power blocks, and the last block is the gas block. The OEF model for each decoupled power block $r$ ( $r = 1,\cdots, N-1$) is presented as follows.

$$\min_{\mathbf{a}_p^r} \sum_{g \in \mathcal{G}_p^r} \mathrm{C}_p(p_g) \tag{1}$$

$$\text{s.t.} \quad \mathrm{P}_g^{\min} \le p_g \le \mathrm{P}_g^{\max}, \quad g \in \mathcal{G}_p^r \tag{2}$$

$$\mathrm{G}_g^{\min} \le g_g \le \mathrm{G}_g^{\max}, \quad g \in \mathcal{G}_g^r \tag{3}$$

$$\theta_i^{\min} \le \theta_i \le \theta_i^{\max}, \quad i \in \mathcal{N}_p^r \tag{4}$$

$$-\mathrm{P}_l \le p_l \le \mathrm{P}_l, \quad l \in \mathcal{L}_p^r \tag{5}$$

$$\mathrm{x}_l \cdot p_l = \begin{cases} \theta_{i(l)} - \theta_{j(l)} & \text{both power nodes are actual} \\ \theta_{i(l)} - \theta_{j'(l)} & \text{one power node is virtual} \end{cases}, \quad l \in \mathcal{L}_p^r \tag{6}$$

$$\sum_{g_p \in \mathcal{G}_p^r} p_{g_p(i)} + \sum_{g_g \in \mathcal{G}_g^r} g_{g_g(i)} + \sum_{l_1 \in \mathcal{L}_p^r} p_{l_1(i)} - \sum_{l_2 \in \mathcal{L}_p^r} p_{l_2(i)} = \sum_{d \in \mathcal{D}_p^r} \mathrm{P}_{d(i)}, \quad i \in \mathcal{N}_p^r. \tag{7}$$

The objective function (1) aims to minimize the generation cost of coal-fired generators, where $\mathrm{C}_p(p_g) = \mathrm{c}_1 \cdot p_g^2 + \mathrm{c}_2 \cdot p_g + \mathrm{c}_3$ is quadratic ($\mathrm{c}_1 > 0$). The generation cost of gas-fired generators is not included. Function (1) is convex (but not necessarily strongly convex).

The feasible region of each power block $r$, $\Omega_p^r$, consists of (2)-(7). Vector $\mathbf{a}_p^r$ is composed of the variables in these constraints. Constraints (2) and (3) enforce the output capacity of coal-fired and gas-fired generators, respectively. Constraint (4) indicates that the phase angle cannot exceed its boundary. Constraint (5) states the thermal limit of a power transmission line. Constraint (6) is The DC power flow equation, where $i(l)$ and $j(l)$ denote a pair of actual power nodes connected by line $l$, while $j(l')$ denotes the virtual power node. This equation is widely employed to characterize the relation between the power flow and the phase angle variables for transmission-level power networks [8]. Constraint (7) is the power balance equation, where $g_p(i)$, $g_g(i)$, $l_1(i)$, $l_2(i)$, and $d(i)$ denote the traditional generator, gas-fired generator, inflow and outflow of power transmission lines, and power loads connecting to power node $i$, respectively. The power balance constraint for virtual power nodes is not included.

### *2.3 Gas Block Model*

The OEF model for the decoupled gas block $r'$ is presented as follows.

$$\min_{\mathbf{a}_g^{r'}} \sum_{w \in \mathcal{W}^{r'}} \mathrm{C}_w \cdot g_w \tag{8}$$

$$\text{s.t.} \quad \mathrm{G}_w^{\min} \le g_w \le \mathrm{G}_w^{\max}, \quad w \in \mathcal{W}^{r'} \tag{9}$$

$$0 \le g_c \le \mathrm{G}_c, \quad c \in \mathcal{C}_g^{r'} \tag{10}$$

$$\mathrm{G}_i^{\min} \le \pi_i \le \mathrm{G}_i^{\max}, \quad i \in \mathcal{N}_g^{r'} \tag{11}$$

$$(g_l)^2 \cdot \mathrm{sgn}(\pi_{i(l)}, \pi_{j(l)}) = \mathrm{W}_l \cdot (\pi_{i(l)} - \pi_{j(l)}), \quad l \in \mathcal{L}_g^{r'} \tag{12}$$

$$\mathrm{sgn}(\pi_{i(l)}, \pi_{j(l)}) = \begin{cases} 1 & \pi_{i(l)} \ge \pi_{j(l)} \\ -1 & \pi_{i(l)} < \pi_{j(l)} \end{cases}, \quad l \in \mathcal{L}_g^{r'} \tag{13}$$

$$\pi_{j(c)} \le \alpha_c \cdot \pi_{i(c)}, \quad c \in \mathcal{C}_g^{r'} \tag{14}$$

$$\sum_{w \in \mathcal{W}^{r'}} g_{w(i)} + \sum_{l_1 \in \mathcal{L}_g^{r'}} g_{l_1(i)} - \sum_{l_2 \in \mathcal{L}_g^{r'}} g_{l_2(i)} + \sum_{c_1 \in \mathcal{C}_g^{r'}} g_{c_1(i)} - \sum_{c_2 \in \mathcal{C}_g^{r'}} g_{c_2(i)} = \sum_{d \in \mathcal{D}_d^{r'}} \mathrm{G}_{d(i)} + \sum_{g'_g \in \mathcal{G}_v^{r'}} \chi_g \cdot g_{g'_g(i)}, \quad i \in \mathcal{N}_g^{r'}. \tag{15}$$

The linear objective function (8) minimizes the output cost of gas wells.

The feasible region for the gas block $r'$, $\Omega_g^{r'}$, consists of (9)-(15). Vector $\mathbf{a}_g^{r'}$ is composed of the variables in these constraints. The output capacity of a gas well is bounded by constraint (9). Constraint (10) restricts the gas flow in the gas compressor (gas active pipeline). Constraint (11) states the limits of the pressure square of the gas node. Equation (12) is called the Weymouth equation, where $i(l)$ and $j(l)$ denote two gas nodes connected by a gas passive pipeline $l$. This (steady-state) equation is derived from its transient-state form, i.e., the momentum equation, and is widely adopted to describe the gas flow in transmission-level gas passive pipelines [34]. Equation (13) depicts its bi-directional property. Constraint (14) is a simplified gas compressor (gas active pipeline) model [35], [36], where $j(c)$ and $i(c)$ represent the outflow and inflow nodes connected by the unidirectional gas compressor $c$. Equation (15) is the gas balance constraint, where $w(i)$, $l_1(i)$, $l_2(i)$, $c_1(i)$, $c_2(i)$, $d(i)$, and $g'_g(i)$ denote the gas well, inflow and outflow of the gas passive pipeline, inflow and outflow of the gas compressor, gas loads, and virtual gas-fired generator connecting to the gas node $i$, respectively.

### *2.4 Coupling Constraints*

In Section 2.1, the IEGS is divided into several blocks by introducing virtual components, i.e., virtual gas-fired generators and virtual power nodes, according to both physical and regional differences. In this section, we formulate the multi-block

IEGS model by exploring the coupling relations between these fully decoupled blocks, and the coupling relations are presented as follows.

$$\theta_{i(l)} = \theta_{i'(l)}, \quad i' \in \bigcup_{r=1}^{N-1} \mathcal{N}_v^r \tag{16}$$

$$g_{g_g(i)} = g_{g'_g(j)}, \quad g'_g \in \mathcal{G}_v^{r'} . \tag{17}$$

Constraint (16) claims that the phase angle of a virtual power node should be identical to that of the corresponding actual power node. For example, assuming that the virtual power node $i'$ in power block $r_2$ is a replica of the actual power node $i$ in another power block $r_1$ (see Fig. 2), the phase angle values of these two nodes should be equal. Similarly, the equivalent relation between virtual and actual gas-fired generators is presented by (17). Notice that the coupling constraints (16)-(17) ensure the equivalence between the original and decoupled (multi-block) models.

Compared with the decoupling method adopted in [19], i.e., removing coupling constraints from the centralized OEF model directly, the proposed method preserves all (corresponding) constraints in each block with the help of virtual components and constructs more simple coupling relations between the decoupled blocks. As a result, the proposed method incorporates fewer system components into coupling constraints and thus requires *less external information* from other blocks. This is significant since this work aims to reduce information exchange from individual agencies and achieve distributed operation for the IEGS.

Note that the values of variables *in the coupling constraints* should not have a significant difference. Otherwise, the variables with large values (e.g., gas-fired power generator outputs) would dominate the small ones (e.g., phase angles), which could result in convergence difficulties. To address this problem, scaling factors are utilized. Specifically, small variables in the coupling constraints are multiplied by scaling factors (i.e., some constants which are larger than 1) to keep their values being of the same order of magnitude or reduce their differences. It is easy to see that the scaling factors are only used in coupling constraints and do not change the variable values, which can be regarded as magnification coefficients.

*2.5 Model Extension*

In practice, we often need to address distributed multi-period OEF problems. The distributed single-period OEF model (1)-(17) can be extended to the multi-period form by adding unit commitment variables, startup and shutdown constraints, and ramping constraints. In addition, the objective functions and constraints in (1)-(17) should be extended to cover multiple periods by adding superscript $t$ (stand for the time period) to corresponding variables. However, the convergence of distributed algorithms remains a challenging issue due to the binary unit commitment variables (the nonconvex Weymouth equation will be convexified in the next section and thus do not affect the algorithm convergence). In order to address this problem, reference [23] adopts an iterative ADMM algorithm. Simulation results indicate that the iterative framework contributes a lot to enhance the convergence, although this algorithm is heuristic and cannot guarantee global optimality. We can use the iterative J-ADMM algorithm (obtained by replacing the ADMM with the J-ADMM in the iterative ADMM algorithm) to solve the distributed multi-period OEF problem.

## 3 Solution Method

This paper adopts the J-ADMM algorithm to solve the distributed OEF problem for the IEGS, and details are presented in Section 3.2. Notice that the convergence of this algorithm is not guaranteed due to nonconvex Weymouth equations in the gas block. In Section 3.1, an ECH-based method is proposed to relax and convexify Weymouth equations. Section 3.3 introduces the solution feasibility check and recovery method.

*3.1 ECH-Based Relaxation Method*

The ECH of a set is convex and contains the convex hull of this set. In addition, an ECH should also: i) contain less redundant elements, and ii) have a simple analytical form.

Denote the convex hull and the ECH of the set $\mathcal{Y}$ as $\mathcal{Y}_{\text{con}}$ and $\mathcal{Y}_{\text{e-con}}$, respectively. According to the above description, $\mathcal{Y}_{\text{e-con}}$ is convex, and the relation $\mathcal{Y}_{\text{con}} \subseteq \mathcal{Y}_{\text{e-con}}$ always holds. The redundant elements refer to the elements which belong to $\mathcal{Y}_{\text{e-con}}$ but do not belong to $\mathcal{Y}_{\text{con}}$. A simple form indicates that we prefer "simple" functions to depict the ECH. For example, a linear function is considered to be simpler than a quadratic function, and a polynomial function is considered to be simpler than exponential, logarithmic, and trigonometric functions. It is easy to see that the proposed ECH is problem-dependent.

The ECH and the convex hull of a Weymouth equation are shown in Fig. 3. The horizontal axis denotes the difference between the pressure squares of gas nodes *i* and *j*. The vertical axis is the gas flow. In this figure, the blue curve depicts the Weymouth equation. In Fig. 3 (a), the region surrounded by red lines is the ECH. It is straightforward to prove that the red lines and the outermost blue lines in Fig. 3 (b) consist of the convex hull boundary.

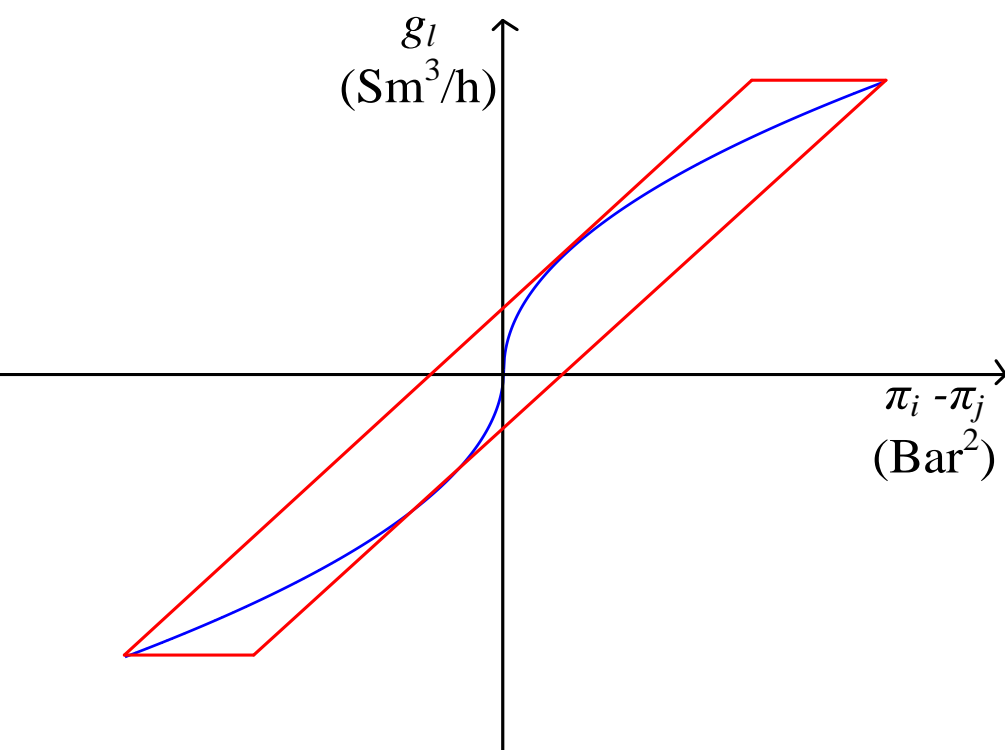

(a) Extended convex hull (ECH) for the Weymouth equation

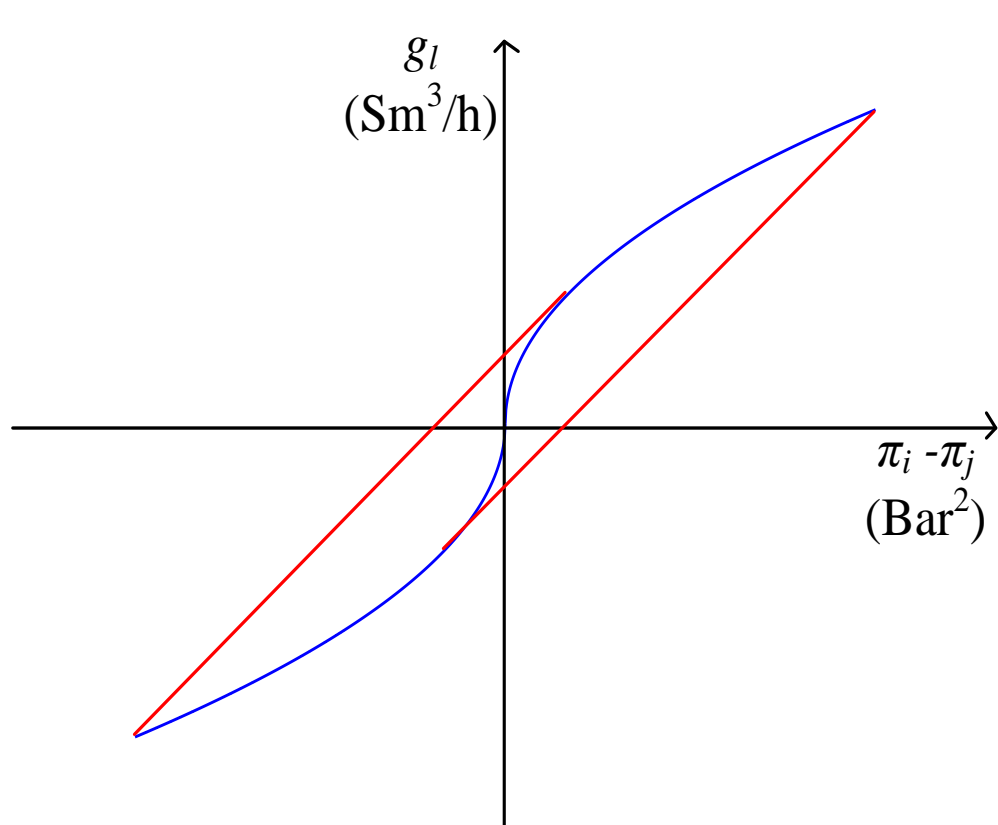

(b) Convex hull for the Weymouth equation

**Fig. 3.** Extended convex hull (ECH) and convex hull for the Weymouth equation

Compared with the tightest convex relaxation, i.e., the convex hull, the proposed ECH in Fig. 3 (a) is more prone to being characterized mathematically, although it contains a larger area consisting of redundant elements. Moreover, the ECH does not introduce any binary variables to the gas block model whilst preserving the bi-directional property of the Weymouth equation. By contrast, the binary variables incurred by the convex hull dissatisfy (one of) the convergence conditions for the distributed algorithm. In this paper, Weymouth equations are replaced by ECH-based constraints, and the mathematical formulation is

$$\mathrm{G}_l^{\min} \le g_l \le \mathrm{G}_l^{\max}, \quad l \in \mathcal{L}_g^{r'} \tag{18a}$$

$$\mathrm{a}_l^{\mathrm{L}} \cdot (\pi_{i(l)} - \pi_{j(l)}) + \mathrm{b}_l^{\mathrm{L}} \le g_l, \quad l \in \mathcal{L}_g^{r'} \tag{18b}$$

$$g_l \le \mathrm{a}_l^{\mathrm{U}} \cdot (\pi_{i(l)} - \pi_{j(l)}) + \mathrm{b}_l^{\mathrm{U}}, \quad l \in \mathcal{L}_g^{r'}. \tag{18c}$$

Constraint (18a) sets the upper and lower bounds of the gas flow in a gas passive pipeline. It also states the upper and lower bounds of the ECH (see Fig. 3 (a)). Constraints (18b) and (18c) represent the left and right bounds of the proposed

ECH (see Fig. 3 (b)), respectively, where $\mathrm{a}_l^{\mathrm{L}}$, $\mathrm{a}_l^{\mathrm{U}}$, $\mathrm{b}_l^{\mathrm{L}}$, and $\mathrm{b}_l^{\mathrm{U}}$ are constants.

Theoretically, besides what has been shown in Fig. 3 (a), the Weymouth equation has the other two forms. These two forms and the corresponding ECHs (being the same as their convex hulls) are presented in Fig. 4 (a) and (b), respectively. The blue line denotes the Weymouth equation. The ECH boundaries consist of blue and red lines. Figure 4 (a) shows the scenario when the upper bound of $\pi_j$ is equal to the lower bound of $\pi_i$, and Figure 4 (b) represents the other scenario when the upper bound of $\pi_j$ is smaller than the lower bound of $\pi_i$. Figure 4 (a) also applies to the scenario when the gas flow direction is fixed. Both ECHs can be uniformly denoted by

$$\mathrm{a}_l^{\mathrm{U}} \cdot (\pi_{i(l)} - \pi_{j(l)}) + \mathrm{b}_l^{\mathrm{U}} \leq g_l, \quad l \in \mathcal{L}_g^{r'} \tag{19a}$$

$$(g_l)^2 \leq \mathrm{W}_l \cdot (\pi_{i(l)} - \pi_{j(l)}), \quad l \in \mathcal{L}_g^{r'}. \tag{19b}$$

Constraints (19a) and (19b) represent the lower and upper boundaries of the proposed ECH, respectively. By replacing the Weymouth equation with ECH-based constraints, the gas block model, which consists of (8)-(11), (14)-(15), and (18) and/or (19), becomes convex.

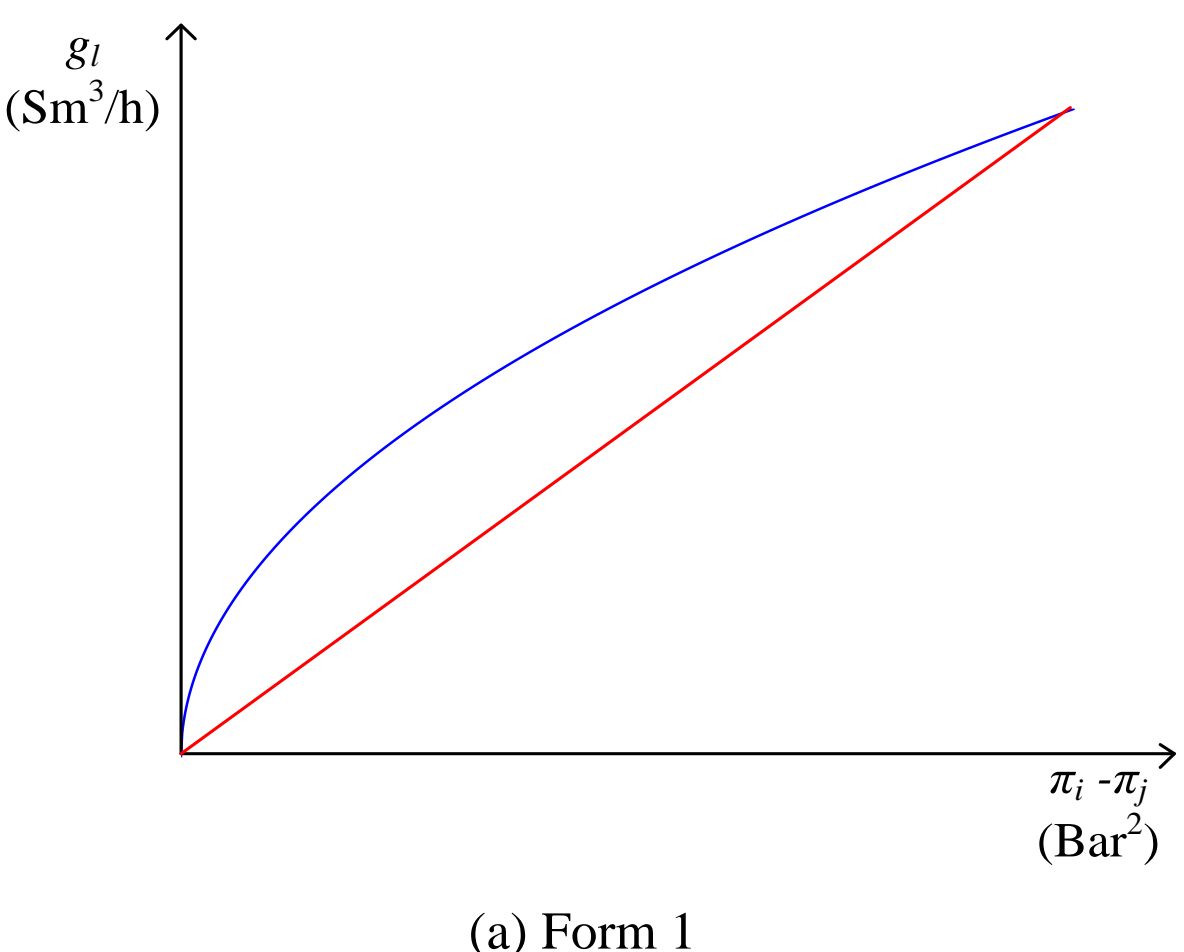

(a) Form 1

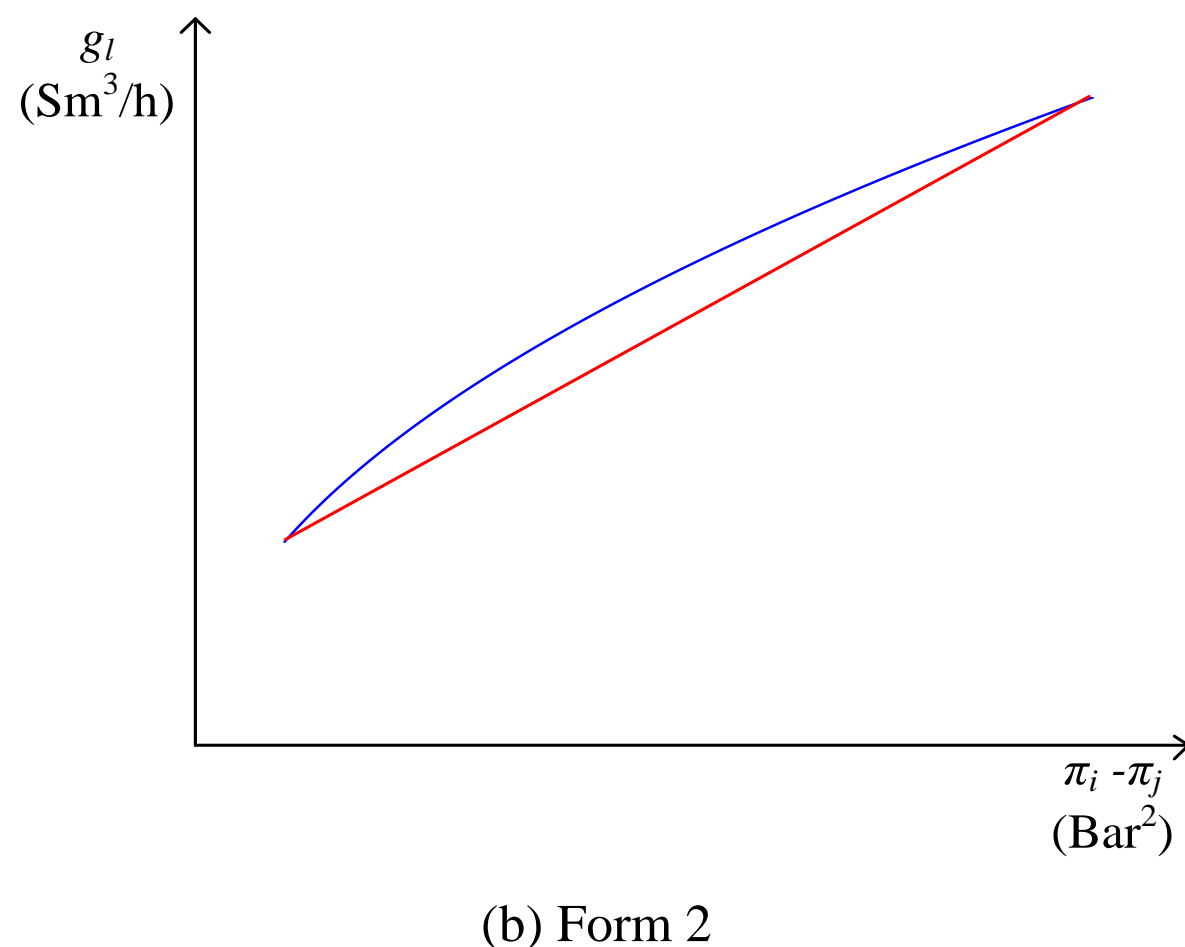

(b) Form 2

**Fig. 4.** Other possible forms of the extended convex hull (ECH) for the Weymouth equation

*3.2 J-ADMM for Multi-Block OEF Problem*

For ease of reading, the compact form of the centralized convexified OEF problem, in which the Weymouth equation is replaced by its ECH-based constraints, is displayed as follows.

$$\min_{\mathbf{x}^1,\cdots,\mathbf{x}^N} f_1(\mathbf{x}^1)+\cdots+f_N(\mathbf{x}^N) \tag{20a}$$

$$\text{s.t.}\ \mathbf{A}_1\cdot\mathbf{x}^1+\cdots+\mathbf{A}_N\cdot\mathbf{x}^N=\mathbf{c} \tag{20b}$$

$$\mathbf{x}^1\in\Omega^1,\cdots,\mathbf{x}^N\in\Omega^N, \tag{20c}$$

where $N$ ($N \geq 2$) is the number of blocks. $\mathbf{x}^r$ are the variables belonging to block $r$ ($r = 1,\cdots, N$). $\Omega^r$ denote convex feasible regions for block $r$. $f_r$ are the convex objective functions for block $r$. $\mathbf{A}_r$ refer to constant matrices, and $\mathbf{c}$ is the constant vector. We do not require that $f_r$ are strongly convex.

The procedures for solving the multi-block optimization problem using the J-ADMM algorithm [28] are as follows.

**Algorithm 1**: J-ADMM algorithm for solving the multi-block distributed OEF problem

1. Initialize penalty parameter d, damping parameter $\gamma$, matrix $\mathbf{P}^r$ ($r = 1,\cdots, N$), stopping criteria $\varepsilon_1$ and $\varepsilon_2$, maximum number of iterations $k_{\max}$, and the number of blocks $N$. Initialize variables $\mathbf{x}_0^r$ ($r = 1,\cdots, N$) and Lagrangian multiplier $\boldsymbol{\lambda}_0$. Set iteration index $k = 0$.

2. Solve (21) for blocks 1 to $N$ in parallel:

$$\mathbf{x}_{k+1}^r=\arg\min_{\mathbf{x}^r\in\Omega_r} f_r(\mathbf{x}^r)+(\mathrm{d}/2)\cdot\left\|\mathbf{A}_r\cdot\mathbf{x}^r+\sum\nolimits_{j\neq r}\mathbf{A}_j\cdot\mathbf{x}_k^j-\mathbf{c}-(\boldsymbol{\lambda}_k/\mathrm{d})\right\|_2^2+(1/2)\cdot\left\|\mathbf{x}^r-\mathbf{x}_k^r\right\|_{\mathbf{P}^r}^2,\quad r=1,\cdots,N, \tag{21}$$

where $\mathbf{x}_k^r$ ($r = 1,\cdots, N$) are the optimal solutions (obtained by (21)) for block $r$ at the ($k$-1)-th iteration. $\|\mathbf{x}_m-\mathbf{x}_n\|_{\mathbf{P}}^2 = (\mathbf{x}_m-\mathbf{x}_n)^{\mathrm{T}}\cdot\mathbf{P}\cdot(\mathbf{x}_m-\mathbf{x}_n)$, and $\mathbf{P}$ is a positive semi-definite matrix. Obtain the updated values $\mathbf{x}_{k+1}^r$ ($r = 1,\cdots, N$) from (21). After obtaining $\mathbf{x}_{k+1}^r$ for $r = 1,\cdots, N$, update Lagrangian multiplier $\boldsymbol{\lambda}_{k+1}$,

$$\boldsymbol{\lambda}_{k+1}=\boldsymbol{\lambda}_k-\gamma\cdot\mathrm{d}\cdot\left(\sum_{r=1}^{N}\mathbf{A}_r\cdot\mathbf{x}_{k+1}^r-\mathbf{c}\right). \tag{22}$$

3. Check whether both the following stop criteria are satisfied:

$$\left\|\sum_{r=1}^{N}\mathbf{A}_r\cdot\mathbf{x}_{k+1}^r-\mathbf{c}\right\|_2^2\leq\varepsilon_1 \tag{23}$$

$$\max\left(\mathrm{d}\cdot\left\|\mathbf{x}_{k+1}^1-\mathbf{x}_k^1\right\|_2^2,\cdots,\mathrm{d}\cdot\left\|\mathbf{x}_{k+1}^N-\mathbf{x}_k^N\right\|_2^2\right)\leq\varepsilon_2. \tag{24}$$

If yes, stop and return $\mathbf{x}_{k+1}^r$ ($r = 1,\cdots, N$); else if $k = k_{\max}$, stop and return NULL (i.e., fail to converge); else, set $k = k + 1$, and go to Step 2.

According to [28], Algorithm 1 converges to its global optimum if $\mathbf{P}^r$ and $\gamma$ satisfy the following conditions.

$$\mathbf{P}^r\succ\mathrm{d}\cdot\left(\frac{1}{\varepsilon_r}-1\right)\cdot\mathbf{E}_r^{\mathrm{T}}\cdot\mathbf{E}_r,\quad r=1,\cdots,N \tag{25}$$

$$\sum_{r=1}^{N}\varepsilon_r<2-\gamma,\quad r=1,\cdots,N. \tag{26}$$

Conditions (25)-(26) can be simplified and transformed into (27) if $\varepsilon_r < (2-\gamma)/N$ ($r = 1,\cdots, N$).

$$\mathbf{P}^r\succ\mathrm{d}\cdot\left(\frac{N}{2-\gamma}-1\right)\cdot\mathbf{E}_r^{\mathrm{T}}\cdot\mathbf{E}_r,\quad r=1,\cdots,N. \tag{27}$$

Compared with other ADMM-based algorithms, the J-ADMM algorithm has the following merits when solving the multi-block OEF problem:

1) *Convergence guarantee*. The global convergence is guaranteed without additional assumptions if the damping parameter $\gamma$ and matrix $\mathbf{P}^r$ ($r = 1,\cdots, N$) are chosen based on (25)-(26) [28]. In contrast, the Gauss-Seidel type ADMM algorithm may not converge when $N \geq 3$ [27], and the Jacobi ADMM algorithm may diverge even when $N = 2$ [28]. Considering the practical situation of the IEGS, i.e., multiple agencies, the number of blocks may be larger than two. The J-ADMM algorithm, which is provably convergent when solving multi-block distributed optimization problems, is quite applicable.

2) *Unique (optimal) operation policy.* The additional proximal term, $(1/2)\cdot||\mathbf{x}^r\text{-}\mathbf{x}_k^r||_{\mathbf{P}^r}^2$, ensures that the mathematical formulation of all blocks is strictly convex. Hence, the optimal solution for the OEF problem is unique. This feature is beneficial for system operators to make (optimal) operation policies for their own blocks, respectively.

3) *Parallel computing*. The J-ADMM algorithm allows the parallel computing (see (21)). In contrast, the Gauss-Seidel type ADMM algorithm can only be implemented in serial. In reality, an IEGS usually consists of multi-blocks, and the J-ADMM algorithm is more efficient to address the distributed OEF problem.

*Remark 1 (parallel computing)*: The J-ADMM allows the parallel computing. Specifically, it updates all blocks simultaneously by solving (21) in parallel, since the J-ADMM algorithm only needs the *k*-th optimal solutions (obtained by solving (21)) of all blocks, i.e., $\mathbf{x}_k^1$, $\cdots$, $\mathbf{x}_k^N$, to update their optimal solutions at the (*k*+1)-th iteration, , i.e., $\mathbf{x}_{k+1}^1$, $\cdots$, $\mathbf{x}_{k+1}^N$,. The parallel computing only works for block updates and does not incorporate the update of the Lagrangian multiplier $\boldsymbol{\lambda}$. In order to elaborate the parallel computing mode of the J-ADMM algorithm, the following flowchart is presented.

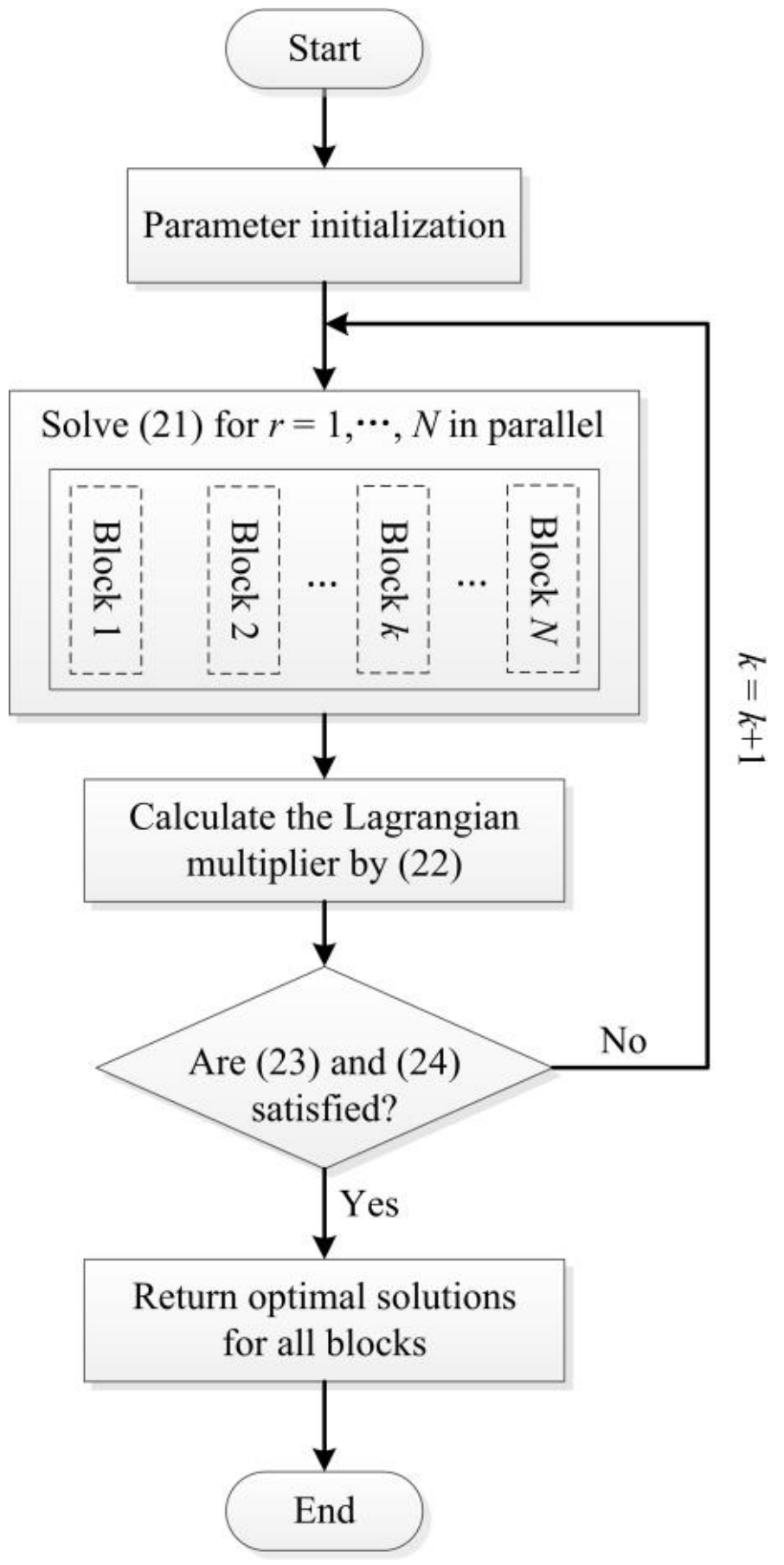

**Fig. 5.** Flowchart of Algorithm 1

As aforementioned, the global convergence of Algorithm 1 is guaranteed if conditions (25) and (26) are satisfied. However, it does not mean that it can converge within a predetermined time. In fact, its computational efficiency is very sensitive to the synergy between the values of the penalty parameter d and damping parameter $\gamma$. Inappropriate combination even leads to the divergence. Specifically, the second and third terms of (21), i.e., $(\mathrm{d}/2)\cdot||\mathbf{A}_r\cdot\mathbf{x}^r + \sum_{j\neq r}\mathbf{A}_r\cdot\mathbf{x}_k^j - \mathbf{c} - (\lambda_k/\mathrm{d})\,||_2^2$ and $(1/2)\cdot||\mathbf{x}^r - \mathbf{x}_k^j||_{\mathbf{P}^r}^2$, correspond to the primal residual in (23) and the dual residual in (24), respectively, and the penalty parameter d and the damping parameter $\gamma$ are their weights ($\mathrm{P}^r$ depends on $\gamma$). The objective function (21) would be dominated by the second term if we chose a large d and a small $\gamma$. Consequently, the primal residual would drop sharply while the dual residual would decrease slowly or even diverge. Similarly, the dual residual would dominate the primal residual with a large $\gamma$ and a small d, which could result in the slow convergence or even divergence of the primal residual. In the case study part (Section 4), we will show how these two parameters affect the algorithm performance. In this work, we use a "grid search" method to choose their values, which is shown as follows.

**Method 1**: Grid search for choosing the values of d and $\gamma$

1. Set two groups of values for the penalty parameter d and the damping parameter $\gamma$. Set $k_{\max}$=1 in Algorithm 1.
2. Choose one pair of untested values for d and $\gamma$ from their groups.
3. Execute Algorithm 1. Record the primal-dual absolute difference (PDAD) value

$$\left| \left\| \sum_{r=1}^{N} \mathbf{A}_r \cdot \mathbf{x}_{k+1}^r - \mathbf{c} \right\|_2^2 \text{-max} \left( \mathrm{d} \cdot \left\| \mathbf{x}_{k+1}^1 - \mathbf{x}_k^1 \right\|_2^2, \cdots, \mathrm{d} \cdot \left\| \mathbf{x}_{k+1}^N - \mathbf{x}_k^N \right\|_2^2 \right) \right|$$

   and the corresponding d and $\gamma$ values, where $|\cdot|$ denotes the absolute value.

4. Check if all pairs of d and $\gamma$ are tested. If yes, go to Step 5; otherwise, go to Step 2.
5. Set d and $\gamma$ to the corresponding values which have the minimum PDAD value.

The idea of Method 1 is to keep the initial primal and dual residuals at the same level to balance their weights (See Section 4.8 for details). Note that Method 1 is a heuristic method to help choose the values for d and $\gamma$ and cannot guarantee any optimum. In addition, the initialization in Step 1 is determined empirically. Recently, some work focuses on the parameter selection for ADMM-based algorithms, such as the adaptive method to update the parameters. Please refer to [16] for details.

*3.3 Solution Feasibility Check and Recovery Method*

The optimal solution obtained by solving the convexified problem may not be feasible for the original nonconvex problem, as the feasible region for the original problem is enlarged. As a result, the optimum of (20) obtained by the distributed optimization algorithm may be smaller than that of the original problem. In order to check whether these two optimums are equal without solving the original nonconvex problem, we have the following proposition.

*Proposition 1 (sufficient condition)*: The original and the convexified problems have the same optimum if the problem (28) is feasible and its objective value is equal to zero.

$$\min_{\boldsymbol{\pi},\, \boldsymbol{\delta}^+,\, \boldsymbol{\delta}^-} \ \mathbf{1}^{\mathrm{T}} \cdot \boldsymbol{\delta}^+ + \mathbf{1}^{\mathrm{T}} \cdot \boldsymbol{\delta}^- \tag{28a}$$

$$\text{s.t.} \ (g_l^*)^2 \cdot \mathrm{sgn}(g_l^*) = \mathrm{W}_l \cdot (\pi_{i(l)} - \pi_{j(l)}), \quad l \in \mathcal{L}_g^{r'} \tag{28b}$$

$$\mathrm{sgn}(g_l^*) = \begin{cases} 1 & g_l^* \geq 0 \\ -1 & g_l^* < 0 \end{cases}, \quad l \in \mathcal{L}_g^{r'} \tag{28c}$$

$$(1-\delta_i^-) \cdot \mathrm{G}_i^{\min} \leq \pi_i \leq (1+\delta_i^+) \cdot \mathrm{G}_i^{\max}, \quad i \in \mathcal{N}_g^{r'} \tag{28d}$$

$$\pi_{j(c)} \leq \alpha_c \pi_{i(c)}, \quad c \in \mathcal{C}_g^{r'} \tag{28e}$$

$$\delta_i^+, \ \delta_i^- \geq 0, \tag{28f}$$

where $\boldsymbol{\delta}^+ = (\delta_1^+, \cdots, \delta_M^+)^{\mathrm{T}}$ and $\boldsymbol{\delta}^- = (\delta_1^-, \cdots, \delta_M^-)^{\mathrm{T}}$ are slack variables. $M = |\mathcal{N}_g^{r'}|$. $g_l^*$ are the gas flow values through gas passive pipelines and are obtained by solving the convexified problem (20).

*Proof*: Assume that $\boldsymbol{\pi}^{**}$ is a feasible solution for (28). Denote the optimal solution for (20) as $\mathbf{x}^*$. $\boldsymbol{\pi}^*$ is the optimal pressure square vector in $\mathbf{x}^*$. Construct a solution, $\mathbf{x}_{opt}$, via substituting $\boldsymbol{\pi}^*$ in $\mathbf{x}^*$ by $\boldsymbol{\pi}^{**}$. $\mathbf{x}_{opt}$ is feasible for the original problem, as it satisfies all the constraints, i.e., (2)-(7) and (9)-(17). $\mathbf{x}_{opt}$ and $\mathbf{x}^*$ are equal except the value of $\boldsymbol{\pi}$. We notice that changing the value of $\boldsymbol{\pi}$ does not influence the objective values of both original and convexified problems. Thus, their optimums are equal. This completes the proof. ■

Proposition 1 also provides a method to recover the feasible optimal solution for the original problem from that for the convexified problem.

*Corollary 1*: Assume that (28) is feasible and its objective value is equal to zero. $\boldsymbol{\pi}^{**}$ is a feasible solution. Denote the optimal solution for (20) as $\mathbf{x}^{*}$. $\boldsymbol{\pi}^{*}$ is the optimal pressure square vector. A solution, $\mathbf{x}_{opt}$, is constructed via replacing $\boldsymbol{\pi}^{*}$ (in $\mathbf{x}^{*}$) by $\boldsymbol{\pi}^{**}$. $\mathbf{x}_{opt}$ is the optimal solution for the original problem.

*Remark 2*: Solving the distributed OEF problem, i.e., the original nonconvex problem, is computationally demanding. Corollary 1 states that if the problem (28), a linear programming problem, is feasible and its objective value is equal to zero, the optimum of the original problem is directly obtained, and its feasible optimal solution can be recovered from the optimal solution for the convexified problem. If the problem (28) is feasible and its objective value is larger than zero, we can quantify the infeasibility degree by observing the values of slack variables $\boldsymbol{\delta}^{+}$ and $\boldsymbol{\delta}^{-}$. In addition, we also obtain a lower bound for the objective value of the original problem.

*Corollary 2*: Problem (28) is always feasible for radial gas networks.

*Proof*: Let each element in $\boldsymbol{\delta}^{+}$ and $\boldsymbol{\delta}^{-}$ equal to a sufficiently large number (smaller than positive infinity). Constraints (28d) and (28f) are always satisfied. For any fixed gas pipeline $l$, $\pi_{i(l)}$ and $\pi_{j(l)}$ need to and only need to satisfy either (28b) or (28e), as a gas pipeline is either a passive pipeline or a compressor (active pipeline). For any fixed gas pipeline $l$ in a radial gas network, we can obtain at least one pair of $\pi_{i(l)}$ and $\pi_{j(l)}$ that satisfies (28b) or (28e) by adjusting the values of $\boldsymbol{\delta}^{+}$ and $\boldsymbol{\delta}^{-}$. By far, all constraints are satisfied (constraint (28c) is a constant). This completes the proof. ■

*Remark 3*: Corollary 2 indicates that for an IEGS with a radial gas network, we can either obtain a feasible optimal solution (for the original problem) or quantify its "infeasibility degree" by the values of slack variables $\boldsymbol{\delta}^{+}$ and $\boldsymbol{\delta}^{-}$. According to [37], gas transmission networks are usually radial, and the proposed solution feasibility check and recovery method has plenty of practical application scenarios.

## 4 Case Study

The proposed methods are tested on two integrated systems, i.e., the 6-node power network with a 7-node gas network and the 118-node power network with a modified 48-node gas network. Detailed system data can be obtained from [38]. Algorithm 1 is coded in MATLAB R2018b. Numerical tests are performed on a PC with an Intel(R) Core(TM) i5-6500 CPU @3.20GHz and a 16 GB memory. The solver used to address the quadratic programming (QP) and the mixed-integer quadratic programming problems is Gurobi 8.1.0. The nonconvex problem (without integer variables) is solved by IPOPT. For the sake of simplicity, condition (27) is adopted in Algorithm 1 to ensure its convergence.

### *4.1 Two-Block Case*

The 6-node power network with a 7-node gas network is divided into two blocks, i.e., one power block and one gas block. The values of the parameters in Algorithm 1 are shown in Table 1, where $\mathbf{M}^{r} = \mathrm{d}\cdot(N/(2-\gamma)-1)\cdot\mathbf{E}_{r}^{\mathrm{T}}\cdot\mathbf{E}_{r}$ ($r = 1,\cdots, N$). In this paper, we set $\mathbf{E}_{r}$ as identity matrices for all $r$. The values of d and $\gamma$ are determined based on Method 1. We do not use scaling factors in this case.

**Table 1.** Parameters setting in Algorithm 1

| Parameter | d | $\gamma$ | $\mathbf{P}^{r}$ | $\varepsilon_1$ | $\varepsilon_2$ | $k_{\max}$ | $N$ |
|---|---|---|---|---|---|---|---|
| Value | 4 | 1 | $1.1\cdot\mathbf{M}^{r}$ | $10^{-4}$ | $10^{-4}$ | 1000 | 2 |

### *4.2 Solution Feasibility Check and Recovery*

In order to validate the effectiveness of the ECH-based relaxation and solution feasibility check and recovery methods, numerical experiments are conducted under different load profiles. Results are shown in Table 2. "Portion" refers to the sum of actual gas well outputs divided by the sum of their maximum output capacities. According to this table, all optimal solutions for the convexified problem are infeasible (short for IF) for the original problem, while feasible optimal solutions (for the original problem) are all successfully recovered (denoted by Y) from infeasible ones, even when gas wells nearly reach their output limits (0.99). These results show the effectiveness of the proposed methods. In the following two-block case studies, power and gas loads are fixed to 286 MW and 6,680 $Sm^3$/h, respectively.

**Table 2.** Recoverability test under different load profiles

| Load profile | Power (MW) | 214 | 243 | 257 | 272 | 286 | 300 | 314 |
|---|---|---|---|---|---|---|---|---|
| | Gas ($Sm^3$/h) | 5010 | 5670 | 6010 | 6340 | 6680 | 7010 | 7350 |
| Portion | | 0.52 | 0.59 | 0.67 | 0.75 | 0.83 | 0.91 | 0.99 |
| Feasibility | | IF | IF | IF | IF | IF | IF | IF |
| Recoverability | | Y | Y | Y | Y | Y | Y | Y |

*4.3 Comparison With Other Methods*

The performances of different methods converting the nonconvex Weymouth equation to enable the distributed operation are compared, and Table 3 shows the results. This including the proposed ECH-based relaxation method, the PWL method with eight divided segments, and the SOC relaxation method. Besides, the unidirectional Weymouth equation, shown as

$$(g_l)^2 = \mathrm{W}_l \cdot (\pi_{i(l)} - \pi_{j(l)}), \quad l \in \mathcal{L}_g^{r'}, \tag{29}$$

is employed to obtain a nonconvex model (referred to as the NCV method). Algorithm 1 is adopted to solve the PWL, SOC, NCV, and ECH (ECH-Jacobi) based models. By contrast, the ECH-based model is also solved by the standard ADMM algorithm with the Gauss-Seidel iteration method (ECH-Gauss) [11]. "F" and "IF" denote feasible and infeasible, respectively. Time refers to the runtime of the entire program.

**Table 3.** Comparison between different methods

| Method | | Optimality ($*10^4$ \$) | Number of Iterations | Time (s) | Feasibility |
|---|---|---|---|---|---|
| PWL | | 7.056 | 31 | 17.0 | F |
| SOC | | 7.056 | 56 | 218.6 | IF |
| NCV | | 7.056 | 30 | 17.5 | F |
| ECH | Jacobi | 7.056 | 30 | 16.8 | F |
| | Gauss | 7.056 | 11 | 9.2 | F |

This test shows that all the methods are feasible, and the same global optimum is obtained. However, the optimal solution for the SOC based model is infeasible due to the second-order cone relaxation of the Weymouth equation [18]. Besides, the computation time of this model is much longer than that of the others. Similar results of the SOC based model can be found in [19]. The ECH based model solved by the Gauss-Seidel iteration method outperforms the others in terms of the running time with the minimum number of iterations, although it only allows serial computing. For this small IEGS, the advantage of Algorithm 1, i.e., parallel computing, is not fully revealed, since there are only two blocks. In Section 4.9, we will conduct another test using a larger system with more blocks.

*4.4 Convergence Performance*

The convergence sequences (primal and dual residuals) of the ECH-Jacobi and ECH-Gauss algorithms are presented in Fig. 6 (a) and (b), respectively. The horizontal axis denotes the number of iterations, and the vertical axis is the residual value. The same sequence is depicted by two kinds of coordinate systems, i.e., Cartesian coordinate system (left) and semi-logarithmic coordinate system (right), to fully capture the convergence trends from distinct perspectives. Compared with the standard ADMM algorithm, the J-ADMM algorithm is more difficult to converge for this case. According to Fig. 6 (a), the primal and dual residuals cannot decrease consistently, although both have a downward trend. In addition, for most of the iterations, they cannot even decrease simultaneously. The possible reason is the alternation of the dominant terms in (21). Specifically, the second and third terms in (21), i.e., $(\mathrm{d}/2)\cdot||\mathbf{A}_r\cdot\mathbf{x}^r + \sum_{j\neq r}\mathbf{A}_r\cdot\mathbf{x}_k^j - \mathbf{c} - (\lambda_k/\mathrm{d})\,||_2^2$ and $(1/2)\cdot||\mathbf{x}^r - \mathbf{x}_k^j||_{\mathbf{P}^r}^2$, are corresponding to the primal residual and the dual residual, respectively. The objective function (21) would be dominated by the primal residual if the dual residual was too small. At this time, the descent of the objective function mainly relied on the primal residual, which would drop sharply. After several iterations, the dual residual would dominate the objective function and affect its decrease. The objective value would continue decreasing by repeating the above procedures, alternatively, which resulted in the inconsistently decrease and convergence difficulty.

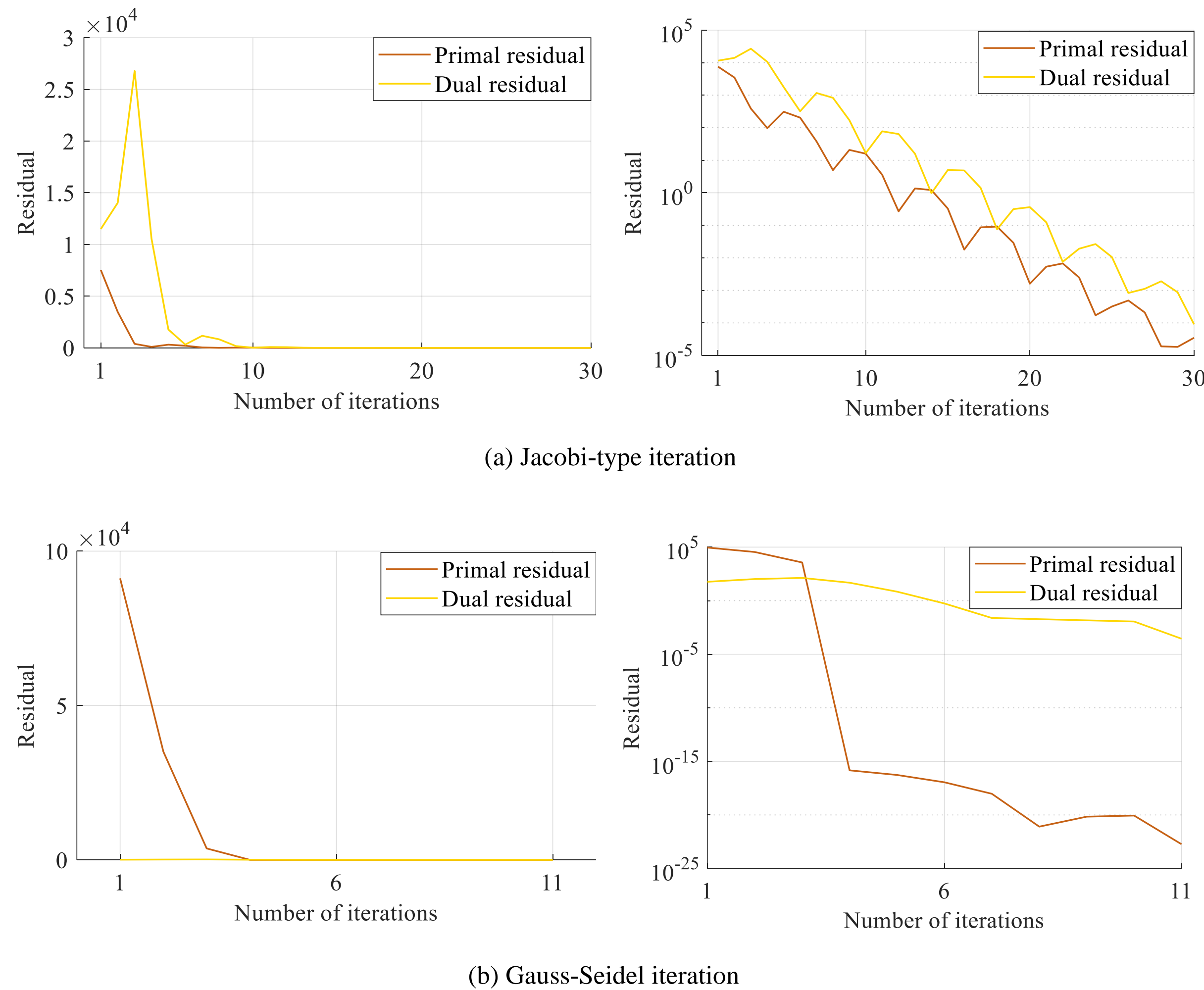

(a) Jacobi-type iteration

(b) Gauss-Seidel iteration

**Fig. 6.** Convergence sequences of different iteration types

### *4.5 Four-Block Case*

In the following subsections, the proposed method is tested on the 118-node power network with a modified 48-node gas network. This IEGS is divided into four blocks, i.e., three power blocks and one gas block. Five compressors are added to gas passive pipelines and lead the gas block to be a "radial" network, as Weymouth equations do not hold for gas active pipelines (please refer to [38] for detailed system topology). We notice that a meshed gas network may have the same property as a radial one due to the existence of compressors. Table 4 shows the initialized parameters in Algorithm 1. The convergence thresholds of both primal and dual residuals are set to $10^{-2}$ to balance the accuracy and the computational cost. The maximum number of iterations is set to $10^4$. The values of d and $\gamma$ are determined based on Method 1. Scaling factors are set to$10^3$.

**Table 4.** Parameters setting in Algorithm 1

| Parameter | d | $\gamma$ | $\mathbf{P}^r$ | $\varepsilon_1$ | $\varepsilon_2$ | $k_{\max}$ | $N$ |
|---|---|---|---|---|---|---|---|
| Value | 1 | 0.2 | $1.1\cdot\mathbf{M}^r$ | $10^{-2}$ | $10^{-2}$ | 1000 | 4 |

*4.6 Solution Feasibility Check and Recovery*

In this subsection, we continue to investigate the effectiveness of the ECH-based relaxation and solution feasibility check and recovery methods. Test results are shown in Table 5. Similar to the results in Section 4.2, feasible optimal solutions for original problems are recovered from infeasible solutions (obtained by the proposed method) under all load profiles. This method still works even for the extreme heavy-load scenario (0.99). We try to analyze why the proposed ECH-based model and the solution recovery method are highly effective for the distributed OEF problem, and the (possibly) sound explanation is as follows. Compared with the underestimation, the overestimation of the transmission capability of the gas passive pipeline is more likely to result in the failure of the proposed solution recovery method (see Fig. 3 (a)). The overestimation means transmitting a specific amount of gas by a smaller pressure difference (compared with the pressure difference calculated by the Weymouth equation (12)-(13) when transmitting the same amount of gas). It is the redundant elements in the ECH that lead to the overestimated transmission capability. In fact, the overestimation caused by constraints (18) is not very large (see Fig. 3 (a)). This overestimation thoroughly vanishes when constraints (19) are adopted (see Fig. 4). Thus, the feasible optimal solutions for the original problem can always be successfully recovered in both cases.

**Table 5.** Recoverability test under different load profiles

| Load profile | Power (MW) | 1250 | 1500 | 1750 | 2125 | 2375 | 2500 |
|---|---|---|---|---|---|---|---|
| | Gas ($Sm^3/h$) | 1350 | 1620 | 1890 | 2295 | 2565 | 2700 |
| Portion | | 0.51 | 0.58 | 0.65 | 0.76 | 0.89 | 0.99 |
| Feasibility | | IF | IF | IF | IF | IF | IF |
| Recoverability | | Y | Y | Y | Y | Y | Y |

*4.7 Gas Transmission Flexibility*

This test aims to observe the gas transmission flexibility in gas passive pipelines, and Table 6 shows the results. Gas flow directions in three pipelines are tracked and pre-defined as A→B, where A and B refer to gas nodes. For example, 11→12 means that the pre-defined gas flow direction is from gas node 11 to gas node 12. Actual gas flow directions under different load profiles that are identical to and opposite to the pre-defined gas flow direction are denoted by "+" and "-", respectively. As shown in this table, gas flow directions vary at different load levels, even when the load level has a slight fluctuation, e.g., pipeline 34→35 with the load increasing from 2375MW and 2565$Sm^3/h$ to 2500MW and 2700$Sm^3/h$. This test shows that ignoring the bi-directional property of the gas flow may lead to a severe discrepancy, which could affect the optimal operation policy for the IEGS. Compared with [19] and [29], the proposed ECH-based method preserves the bi-directional property of gas flows in gas passive pipelines and leads to a more flexible and effective gas block model. In the following four-block case studies, power and gas loads are fixed to 2500 MW and 1580 $Sm^3/h$, respectively.

**Table 6.** Gas flow direction under different load profiles

| Load profile | Power (MW) | 1250 | 1500 | 1750 | 2125 | 2375 | 2500 |
|---|---|---|---|---|---|---|---|
| | Gas ($Sm^3$/h) | 1350 | 1620 | 1890 | 2295 | 2565 | 2700 |
| Transmission direction of gas pipeline | 11→12 | + | + | + | - | + | + |
| | 34→35 | + | + | + | + | + | - |
| | 37→38 | - | - | + | - | + | + |

*4.8 Computation Time*

As aforementioned, the values of the penalty parameter d and the damping parameter $\gamma$ greatly influence the computational efficiency of Algorithm 1. Theoretically, a big penalty term, i.e., $(d/2)\cdot||\mathbf{A}_r\cdot\mathbf{x}^r + \sum_{j\neq r}\mathbf{A}_r\cdot\mathbf{x}_k^j - \mathbf{c} - (\lambda_k/d)\,||_2^2$, would be added to the objective function (21) if the value of d were too large. As a result, the dual residual term, i.e. $(1/2)\cdot||\mathbf{x}^r - \mathbf{x}_k^j||_{\mathbf{P}^r}^2$ in (21), would decrease slowly. However, if the value of d were too small, the decline rate of the primal residual term, i.e., $(d/2)\cdot||\mathbf{A}_r\cdot\mathbf{x}^r + \sum_{j\neq r}\mathbf{A}_r\cdot\mathbf{x}_k^j - \mathbf{c} - (\lambda_k/d)\,||_2^2$, would become very slow. Thus, this penalty parameter needs to be chosen within a rational range. In addition, according to (27) and Table 4, the values of matrices $\mathbf{P}^r$ ($r = 1,\cdots, N$) are decided by $\gamma$ and closely related to the dual residual term $(1/2)\cdot||\mathbf{x}^r - \mathbf{x}_k^j||_{\mathbf{P}^r}^2$ in the objective function (21). Thus, the selection of $\gamma$ also impacts algorithm efficiency. Simulation results obtained by solving the ECH-based model using Algorithm 1 are shown in Table 7. Computation time refers to the runtime of the entire program. Test results indicate that the penalty parameter d and the damping parameter $\gamma$ play a critical role in deciding the running time of the proposed algorithm. In addition, they are quite sensitive. Namely, slight changes in their values may lead to dramatic fluctuations in the computation time. In order to improve the computational efficiency of Algorithm 1, both parameters should be designed delicately.

**Table 7.** Influence of parameters on computation time (s)

| $\gamma$ \ d | 0.1 | 0.5 | 1 | 2 | 10 |
|---|---|---|---|---|---|
| 0.01 | >3600.0 | 2812.0 | 2751.8 | 2217.6 | >3600.0 |
| 0.1 | 1597.3 | 934.2 | 990.1 | 1114.6 | >3600.0 |
| 0.2 | 1236.4 | 618.1 | 507.3 | 895.0 | 2289.0 |
| 1.0 | 1410.4 | 805.3 | 1389.3 | 1407.2 | >3600.0 |
| 1.5 | 959.4 | 1604.2 | 1387.1 | 1849.4 | >3600.0 |
| 1.9 | >3600.0 | >3600.0 | >3600.0 | >3600.0 | >3600.0 |

*4.9 Comparison With Other Methods*

Similar to the comparison in Section 4.3, different model converting methods, which enable the distributed operation, are tested using this larger IEGS, and results are shown in Table 8. Algorithm 1 is used to solve PWL, SOC, NVC, and ECH based models. Besides, the direct extension of the standard ADMM algorithm is used to solve the proposed multi-block problem (written as ECH-Gauss). Time refers to the runtime of the entire program. Binary variables are introduced to PWL and SOC based models, leading to nonconvex mixed-integer programming problems, and both fail to converge within ten hours. The NCV-based model does not converge within ten hours due to nonconvex Weymouth equations. The proposed ECH-based model (leading to a convex QP problem) converges under both algorithms. Algorithm 1, which allows parallel computing, outperforms the ECH-Gauss algorithm on the running time, although the former takes more iterations to converge. This test indicates that Algorithm 1 has great potential to address large-scale multi-block distributed optimization problems, especially when combining with state-of-the-art parallel computing techniques.

**Table 8.** Comparison between different methods

| Method | | Optimality (*$10^4$ $) | Number of Iterations | Time (s) | Feasibility |
|---|---|---|---|---|---|
| PWL | | - | - | >36000.0 | - |
| SOC | | - | - | >36000.0 | - |
| NCV | | - | - | >36000.0 | - |
| ECH | Jacobi | 13.400 | 976 | 507.3 | F |
| | Gauss | 13.400 | 789 | 1541.1 | F |

## 5 Conclusion

This paper proposes an extended convex hull-based multi-block electricity-gas system model to address the optimal energy flow problem in a distributed manner. The nonconvex Weymouth equation is convexified and replaced by the extended convex hull-based constraints. The Jacobi-Proximal alternating direction method of multipliers algorithm is adopted to solve this convexified problem. After that, the feasibility of the optimal solution for the convexified problem is checked, and a sufficient condition is developed. The optimal solution for the original nonconvex problem is recovered from that for the convexified problem if the sufficient condition is satisfied. Simulation results show that:

1) The proposed extended convex hull-based model preserves the bi-directional property of the Weymouth equation and successfully captures the changes of gas flow directions, even under slight load fluctuations. In addition, this model shows its superiority of the convergence and computation time over other models.
2) The Jacobi-Proximal alternating direction method of multipliers algorithm is effective and efficient in solving multi-block distributed optimal energy flow problems. Notice that the penalty parameter d and the damping parameter $\gamma$ should be designed delicately to further improve computational efficiency.
3) Though only a sufficient condition, the solution recovery method is effective in helping obtain a feasible optimal solution for radial gas networks under different load profiles.
4) The proposed framework can address distributed optimal energy flow problems with multiple blocks and thus has huge potential to applies to distributed optimization problems for realistic large-scale multi-block integrated electricity and gas systems, especially combining with the state-of-the-art parallel computing techniques.

Future work includes developing a new solution recovery method for meshed gas networks.

## 6 Acknowledgment

This work was supported in part by the National Natural Science Foundation of China under Grant 51677160, in part by the Research Grants Council of Hong Kong under Grant GRF17207818, in part by the Research Grants Council of Hong Kong through the Theme-based Research Scheme under Project No. T23-701/14-N, in part by the U.S. National Science Foundation under Grant ECCS-1552073, and in part by the U.S. Department of Energy under Award DE-EE0007998.

## 7 References

[1] Middleton RS, Gupta R, Hyman JD, Viswanathan HS. The shale gas revolution: Barriers, sustainability, and emerging opportunities. Appl Energy 2017;199(1):88-95.

[2] U.S. Energy Information Administration. Annual Energy Review, 2019 [Online]. Available: https://www.eia.gov/totalenergy/data/annual/index.Php

[3] National Statistics. Electricity Generation, Trade and Consumption, 2019. [Online]. Available: https://www.gov.uk/government/statistics/electricity-section-5-energy-trends

[4] Qiu J, Dong ZY, Zhao JH, Meng K, Zheng Y, Hill DJ. Low carbon oriented expansion planning of integrated gas and power systems. IEEE Trans Power Syst 2015;30(2):1035-1046.

[5] Zeng Q, Zhang B, Fang J, Chen Z. A bi-level programming for multistage co-expansion planning of the integrated gas and electricity system. Appl Energy 2017;200:192-203.

[6] Chen S, Wei Z, Sun G, Sun Y, Zang H, Zhu Y. Optimal power and gas flow with a limited number of control actions. IEEE Trans Smart Grid 2018;9(5):5371-5380.

[7] Belderbos A, Bruninx K, Delarue E, D'haeseleer W. Facilitating renewables and power-to-gas via integrated electrical power-gas system scheduling. Appl Energy 2020;275:54-68.

[8] Correa-Posada CM, Sánchez-Martin P. Security-constrained optimal power and natural-gas flow. IEEE Trans Power Syst 2014;29(4):1780-1787.

[9] He C, Dai C, Wu L, Liu T. Robust network hardening strategy for enhancing resilience of integrated electricity and natural gas distribution systems against natural disasters. IEEE Trans Power Syst 2018;33(5):5787-5798.

[10] Farrokhifar M, Nie Y, Pozo D. Energy systems planning: A survey on models for integrated power and natural gas networks coordination. Appl Energy 2020;262:114567.

[11] Boyd S, Parikh N, Peleato B, Eckstein J. Distributed optimization and statistical learning via the alternating direction method of multipliers. Found Trends in Mach Learn 2011;3(1):1-122.

[12] Erseghe T. Distributed optimal power flow using ADMM. IEEE Trans Power Syst 2014;29(5): 2370-2380.

[13] Fang X, Hodge BM, Jiang H, Zhang Y. Decentralized wind uncertainty management: Alternating direction method of multipliers based distributionally-robust chance constrained optimal power flow. Appl Energy 2019;239: 938-947.

[14] Khaki B, Chu C, Gadh R. Hierarchical distributed framework for EV charging scheduling using exchange problem. Appl Energy 2019;241:461-471.

[15] Lai K, Illindala MS. A distributed energy management strategy for resilient shipboard power system. Appl Energy 2018;228:821-832.

[16] Mavromatis C, Foti M, Vavalis M. Auto-tuned weighted-penalty parameter ADMM for distributed Optimal Power Flow. IEEE Trans Power Syst (early access).

[17] Wen Y, Qu X, Li W, Liu X, Ye X. Synergistic operation of electricity and natural gas networks via ADMM. IEEE Trans Smart Grid 2018;9(5):4555-4565.

[18] Borraz-Sanchez C, Bent R, Backhaus S, Hijazi H, Hentenryck P. Convex relaxations for gas expansion planning. INFORMS J Comput 2016;28(4):645-656.

[19] Wang C, Wei W, Wang J, Bai L, Liang Y, Bi T. Convex optimization based distributed optimal gas-power flow calculation. IEEE Trans Sustain Energy 2018;9(3):1145-1156.

[20] Lipp T, Boyd S. Variations and extension of the convex-concave procedure. Optim Eng 2016;17:263-287.

[21] He C, Wu L, Liu T, Shahidehpour M. Robust co-optimization scheduling of electricity and natural gas systems via ADMM. IEEE Trans Sustain Energy 2017;8(2):658-670.

[22] Correa-Posada CM, Sanchez-Martin P. Gas network optimization: A comparison of piecewise linear models. 2014 [Online]. Available: http://www.optimization-online.org/DB_HTML/2014/10/4580.html

[23] Takapoui R, Moehle N, Boyd S, Bemporad A. A simple effective heuristic for embedded mixed-integer quadratic programming. Am Control Conf (ACC) 2016. p. 5619-5625.

[24] He Y, Yan M, Shahidehpour M, Li Z, Guo C, Wu L, et al. Decentralized optimization of multi-area electricity-natural gas flows based on cone reformulation. IEEE Trans Power Syst 2018;33(4):4531-4542.

[25] Qu K, Shi S, Yu T, Wang W. A convex decentralized optimization for environmental-economic power and gas system considering diversified emission control. Appl Energy 2019;240:630-645.

[26] Qi F, Shahidehpour M, Li Z, Wen F, Shao C. A chance-constrained decentralized operation of multi-area integrated electricity-natural gas systems with variable wind and solar energy. IEEE Trans Sustain Energy 2020;11(4): 2230-2240.

[27] Chen C, He B, Ye Y, Yuan X. The direct extension of ADMM for multi-block convex minimization problems is not necessarily convergent. Math Program 2014;155(1):57-79.

[28] Deng W, Lai MJ, Peng Z, Yin W. Parallel multi-block ADMM with o(1/k) convergence. J Sci Comput 2017;71(2):712-736.

[29] Chen S, Conejo AJ, Sioshansi R, Wei Z. Unit commitment with an enhanced natural gas-flow model. IEEE Trans Power Syst 2019;34(5):3729-3738.

[30] Kou X, Li F, Dong J, Olama M, Starke M, Chen Y, et al. A comprehensive scheduling framework using SP-ADMM for residential demand response with weather and consumer uncertainties. IEEE Trans Power Syst (early access).

[31] Chen Z, Li Z, Guo C, Wang J, Ding Y. Fully distributed robust reserve scheduling for coupled transmission and distribution systems. IEEE Trans Power Syst (early access).

[32] Lavaei J, Low SH. Zero duality gap in optimal power flow problem. IEEE Trans Power Syst 2012;27(1):92-107.

[33] Farivar M, Low SH. Branch flow model: relaxations and convexification—Part I. IEEE Trans Power Syst 2013;28(3):2554-2564.

[34] LMNO Engineering. Weymouth, panhandle A and B equations for compressible gas flow. [Online]. Available: https://www.lmnoeng.com/ Flow/weymouth.php

[35] Wang C, Wei W, Wang J, Bi T. Convex optimization based adjustable robust dispatch for integrated electric-gas systems considering gas delivery priority. Appl Energy 2019;239:70-82.

[36] Wolf D, Smeers Y. The gas transmission problem solved by an extension of the simplex algorithm. Manage Sci 2000;46(11):1454–1465.

[37] Wang C, Gao R, Wei W, Shafie-khah M, Bi T, Catalão JPS. Risk-based distributionally robust optimal gas-power flow with Wasserstein distance. IEEE Trans Power Syst 2019;34(3):2190-2204.

[38] [Online]. Available: https://sites.google.com/site/rongpengliu1991/